\documentclass[12pt,reqno]{amsart}
\usepackage{amssymb,amsmath,amsthm}

\newtheorem{theorem}{Theorem}

\newtheorem{corollary}{Corollary}

\theoremstyle{definition}

\def\eps{\epsilon}

\DeclareMathOperator{\Rre}{Re}
\DeclareMathOperator{\Imm}{Im}
\def\bar{\overline}
\def\w{\widehat}

\def\eps{\epsilon}
\def\sgn{\textrm{sgn}}
\newcommand{\R}{\mathbb{R}}
\newcommand{\C}{\mathbb{C}}
\newcommand{\eplj}{\epsilon^\lambda_j}
\newcommand{\mulj}{\mu^\lambda_j}
\newcommand{\epj}{\epsilon_j^\lambda}
\newcommand{\atopp}[2]{\genfrac{}{}{0pt}{2}{#1}{#2}}
\newcommand{\I}{\mathcal{I}}
\newcommand{\dbar}{\bar\partial}
\newcommand{\p}{\partial}
\newcommand{\z}{\bar z}

\begin{document}

\title{The $\Box_b$-heat equation on quadric manifolds}
\author{Albert Boggess and Andrew Raich}

\thanks{The second author is partially funded by NSF grant DMS-0855822}

\address{Department of Mathematics \\ Texas A\&M University\\ Mailstop 3368 \\ College Station, TX  77845-3368 
\newline\newline
Department of  Mathematics \\ 1 University of Arkansas \\ SCEN 327 \\ Fayetteville, AR 72701}


\email{boggess@tamu.edu, araich@uark.edu}

\begin{abstract}In this article, we give an explicit calculation of the partial Fourier transform of the $\Box_b$-heat equation on quadric submanifolds
of $M\subset\C^n\times\C^m$. As a consequence, we can also compute the heat kernel associated to the weighted $\dbar$-equation in $\C^n$ when the weight 
is given by $\exp(-\phi(z,z)\cdot\lambda)$ where
$\phi:\C^n\times\C^n\to \C^m$ is a quadratic, sesquilinear form and $\lambda\in\R^m$. Our method involves the representation theory of the Lie group $M$ and the
group Fourier transform.
\end{abstract}

\maketitle

%%%%%%%%%%%%%%%%%%%%
%
%	INTRODUCTION
%
%%%%%%%%%%%%%%%%%%%%%
\section{Introduction}\label{sec:intro}

The purpose of this article is to present an explicit calculation of the Fourier transform of the fundamental solution of the $\Box_b$-heat equation
on quadric submanifolds $M\subset \C^n\times\C^m$. A quadric submanifold can be thought of as a generalization of the Heisenberg group -- it is a 
Lie group with a known representation theory \cite{PeRi03}, and the technique of using Hermite functions to compute the heat kernel, as done in
\cite{Hul76, BoRa09} and elsewhere, can be extended to work in this situation as well. 

A consequence of our fundamental solution computation is that
we can explicitly compute the heat kernel associated to weighted $\dbar$-problem in $\C^n$ when the weight is given by $\exp(-\phi(z,z)\cdot\lambda)$ where
$\phi:\C^n\times\C^n\to \C^m$ is a quadratic, sesquilinear form and $\lambda\in\R^m$. 
This computation partially generalizes the results in  \cite{BoRa09}.
When $m=1$ and the weight is given by $\exp(\tau P(z_1, \dots, z_n))$ where $\tau\in\R$, 
$P(z_1,\dots,z_n) = \sum_{j=1}^n p_j(z_j)$, and  $p_j$ are subharmonic, nonharmonic polynomials, 
Raich \cite{Rai06h, Rai06f, Rai07, Rai09h} has estimated the heat kernel associated to the weighted $\dbar$-problem. If, in addition, $n=1$, 
the weighted $\dbar$-problem and explicit construction of Bergman and Szeg\"o kernels have been studied by a number of authors in different contexts, e.g.,
\cite{Christ91,Has94,Has95, Has98,FoSi91, Ber92}.
We also note that quadric manifolds are related to  $H$-type groups on which Yang and Zhu  have computed the heat kernel for the sub-Laplacian \cite{YaZh08}.

The remainder of the paper is organized as follows: in Section \ref{sec:definitions, results}, we define our terms and state our main results. Section
\ref{sec:rep theory} provides the necessary background from representation theory. In
Sections \ref{sec:heat equation and Group F.T.} and \ref{sec:heat kernel computation}, we apply the representation theory to the heat kernels and prove the main results.

%%%%%%%%%%%%%%%%%%%%%
%
%	Quadric Submanifolds and the $\Box_b$-heat equation
%
%%%%%%%%%%%%%%%%%%%%%%
\section{Quadric Submanifolds and the $\Box_b$-heat equation.} \label{sec:definitions, results}

%
%	quadric submanifolds
%
\subsection{Quadric submanifolds}
Let $M$ be the the quadric submanifold in
$\C^n \times \C^m$ defined by
\[
M= \{(z,w) \in \C^n \times \C^m; \ \Imm { w} = \phi(z,z) \}
\]
where $\phi: \C^n \times \C^n \mapsto \C^m $ is a sesquilinear
form (i.e. $\phi (z, z') = \overline{\phi (z',z)}$). For emphasis, we sometimes write $M_\phi$ to denote the dependence
of $M$ on the quadratic function $\phi$. Note that $M_{-\phi}$
is biholomorphic to $M_\phi$ by the change of variables
$(z,w) \mapsto (z,-w)$.

For 
$\lambda \in \R^m$, let
\[
\phi^\lambda (z,z') = \phi (z,z') \cdot \lambda
\]
where $\cdot$ is the ordinary dot product (without conjugation).
Observe that $\phi^\lambda (z,z') $
is a sesquilinear scalar-valued form with an associated Hermitian
matrix. Let $v^\lambda_1, \dots, v^\lambda_n$
be an orthonormal basis for $\C^n$ with 
\[
\phi^\lambda (v^\lambda_j, v^\lambda_k) = \delta_{jk} \mu_j (\lambda)
\]
where $\mu_j (\lambda) = \mulj$ are the eigenvalues of 
the matrix associated with $\phi^\lambda$.

%
%	the Lie Group structure
%
\subsection{Lie Group Structure.}
By projecting $M \subset \C^n \times \C^m$ onto $G=\C^n \times \R^m$,
the Lie group structure of $M$ is isomorphic to the following
group structure on $G$:
\[
g g'=(z,t)  (z', t') =\big(z+z', t+t' + 2\Imm \phi (z,z') \big).
\]
Note that $(0,0)$ is the identity in this group
structure and that the inverse of $(z,t)$ is $(-z,-t)$.

The \emph{right invariant} vector fields are 
given as follows: let $g \in G$; if $X$ is a vector field,
then we denote its value at $g$ by $X(g)$ as a element
of the tangent space of $M$ at $g$. 
Define
$R_g: G \mapsto G$ by $R_g(g') = g'g$;
then the right invariant vector fields, $X(g)$, are obtained
by pushing forward the vectors in the tangent space at the origin via
the differential of the map $R_g$. In particular,
a vector field $X$ is right invariant if and only if
$X(g)=(R_g)_* \{ X(0) \}$, where $(R_g)_* $ denotes
the push forward operator. Let 
$v $ be a vector in $\C^n\approx \R^{2n}$
which can be identified with the tangent space
of $M$ at the origin. Let $\partial_v$ be the real
vector field given by the directional 
derivative in the direction of $v$. Then
the right invariant vector field at an arbitrary
$g=(z,w) \in M$ corresponding to $v$ is given by
\[
X_v (g)=  \partial_v + 2 \Imm \phi (v,z) \cdot D_t
=\partial_v - 2 \Imm \phi (z,v) \cdot D_t
\]
where $D_t = (\p_{t_1},\dots, \p_{t_m})$,
(see Section 1 in Peloso/Ricci  \cite{PeRi03}). 
Let $Jv$ be the vector in $\R^{2n}$ which corresponds
to $iv$ in $\C^n$ (where $i=\sqrt{-1}$).
The CR structure on $G$ is then spanned by vectors of the 
form:
\[
Z_v(g)=(1/2)(X_v-iX_{Jv})
=(1/2)(\partial_v - i \partial_{Jv}) -i \bar{\phi (z,v)} \cdot D_t 
\]
and
\[
\bar Z_v(g)=(1/2)(X_v+iX_{Jv})
=(1/2)(\partial_v + i \partial_{Jv}) +i \phi (z,v) \cdot D_t. 
\]
Also,
\[
[X_v, X_{v'} ] = 4 \Imm \phi (v',v) \cdot D_t, \ \ 
[Z_v, Z_{v'}] = 0 
\] 
and 
\[
[\bar Z_v, \bar Z_{v'}] =0, \ \ \  
[Z_v, \bar Z_{v'}] = 2i \phi (v, v') \cdot D_t.
\]
We often drop the $g$ in the vector field notation.
The vector field definition of the Levi form of $M$
is the map $v \mapsto \textrm{proj}( [Z_v, \bar Z_v])$, where 
proj stands for the projection onto the totally real
part of the tangent space of $M$ at the origin (i.e. the $t$-axis).
From the above equation, clearly the Levi form of $M$ can be identified
with the map $v \mapsto \phi(v,v)$, as mentioned at the beginning
of this section. 

Recall that for any $\lambda \in \R^m$, the set of vectors
$v^\lambda_1, \dots, v^\lambda_n$ is an orthonormal basis
which diagonalizes $\phi^\lambda (z,z) = \phi(z,z) \cdot \lambda$.
For $\lambda\in\R^m$, define the function ${\nu(\lambda)}$ by
\[
{\nu(\lambda)} = \text{rank}(\phi^\lambda).
\]
The function ${\nu(\lambda)}$ satisfies $0\leq {\nu(\lambda)} \leq n$ and as in \cite{PeRi03}, 
\[
\{ \lambda\in \R^m :  \nu(\lambda) \equiv \max_{\tilde \lambda\in\R^m} \nu(\tilde\lambda) \}
\]
is a Zariski-open set $\Omega\subset\R^m$ that carries full measure, i.e., $|\R^m \setminus \Omega|=0$. 
We identify $x$ with $(x^\lambda_1, \dots, x^\lambda_n)$
and $y$ with $(y^{\lambda}_1, \dots, y^{\lambda}_n)$
We also write $z=\sum_{j=1}^n (x^\lambda_j +i y^\lambda_j) v^\lambda_j$
for $z = x+iy \in \C^n$. Additionally, we let $z' = (z_1^\lambda,\dots, z_{\nu(\lambda)}^\lambda)$, $z'' = (z_{\nu(\lambda)+1}^\lambda,\dots,z_n^\lambda)$ 
and similarly for $x$ and $y$.

Since the right invariant vector fields corresponding
to $\phi$ are equal to the left invariant vector fields
corresponding to $-\phi$ and $M_{-\phi}$
is biholomorphic to $M_\phi$, any analysis involving 
right invariant vector fields yields corresponding information
about the left invariant vector fields and vice versa.

%
%	COMPUTATION OF \Box_b
%
\subsection{$\Box_b$ Calculations}
Let $v_1, \dots  , v_n$ be any orthonormal basis for $\C^n$.
Let $X_j= X_{v_j}$, $Y_j=X_{Jv_j}$,
and let $Z_j=(1/2) (X_j-iY_j)$, $\bar{Z}_j = (1/2) (X_j+iY_j)$
be the right invariant vector fields defined above 
(which are also the left invariant vector fields for the group structure
with $\phi$ replaced by $-\phi$).
Also let $dz_j$ and $d\bar{z}_j$ be the dual basis.
A $(0,q)$-form can be expressed as $\sum_{K\in \I_q} \phi_K\, d\z^K $
where $\I_q = \{ K= (k_1,\dots,k_q) : 1 \leq k_1 < \cdots < k_q \leq n\}$.
Proposition 2.1 in  \cite{PeRi03} states that 
\[
\Box_b ( \sum_{K\in\I_q} \phi_K\, d\z^K )
= \sum_{K,L\in I_q}  \Box_{LK} \phi_K \, d \bar{z}^L
\]
where
\begin{equation}
\label{boxLK}
\Box_{LK} = - \delta_{LK} \mathcal{L} + M_{LK}
\end{equation}
with the sub-Laplacian on $G$
\[
\mathcal{L} = (1/2) \sum_{k=1}^n \bar Z_k Z_k + Z_k \bar Z_k
\]
and 
\[
M_{LK} = \left\{
\begin{array}{cc} \displaystyle
\vspace{.05in} \frac 12 \left( \sum_{k \in K} [Z_k , \bar Z_k] - \sum_{k \not\in K} [Z_k , \bar Z_k] 
\right) & \textrm{if} \ K=L \\
\vspace{.1in} \eps(K,L) [Z_k, \bar Z_l] & \textrm{if} \ |K \cap L | = q-1 \\
0 & \textrm{otherwise}
\end{array}
\right.
\]
Here, $\eps(K,L) $ is $(-1)^{d}$ where 
$d$ is the number of elements in $K \cap L$ between
the unique element $k \in K-L$ and the unique element
$l \in L-K$.
The above theorem is stated and proved in  \cite{PeRi03} for the 
left-invariant vector fields. If right invariant
vector fields are used, then the above theorem provides a formula for 
$\Box_b$ associated to $M_{-\phi}$.

%From direct calculation in the case of Example I (the quadric in $\C^4$ with 
%surjective Levi form), we have $\sum_{k=1}^2 [Z_k, \bar Z_k] =0$,
%for $Z_k=Z_{v_k}$ where $\lambda = (1,0)$, $v_1=(1,0)$, and $v_2=(0,1)$. So the above
%formula for  
%$\Box_b$
%is just the sub-Laplacian
%in the case when $q=0$ or 2 (the degrees where $\Box_b$ 
%is always solvable).

For later, we record the diagonal part of $\Box_b$, i.e., 
$\Box_{LL}$. Using (\ref{boxLK}) with $L=K$ and the above formulas for
$Z_k$, we obtain
\begin{multline}
\label{eqn:BoxKK}
\Box_{LL} = -\frac 14 \Delta
+2 \Imm \left\{ \sum_{k=1}^n \phi(z,v_k)\p_{z_k} \right\} \cdot D_{t}
- \sum_{k=1}^n \big( \phi(z,v_k) \cdot D_t \big) \big( \bar{\phi(z,v_k)} \cdot D_t \big) \\
+i \left( \sum_{k \in L} \phi(v_k, v_k) \cdot D_t 
- \sum_{k \not\in L} \phi(v_k, v_k) \cdot D_t \right)
\end{multline}
where $\Delta$ is the usual Laplacian in the $z$-coordinates.
% and 
%$D_{tz_k}$
%is the vector of differential operators whose $j$th component is 
%$\partial^2/\partial_{t_j} \partial_{z_k}$.
For example, in the classic case of  the Heisenberg group,
$\phi(z,z)=|z|^2$, and $Z_k=\p_{z_k}-i \bar{z_k} \p_t$, and
$\Box_b$ is a diagonal operator (since $[Z_k, \bar{Z}_l] =0$
when $k \not= l$). The above formula for $\Box_{LL}$ then
gives the coefficient of $\Box_b$ acting on forms of the type
$\phi_L (z) d \bar{z}^L$.

%
%	the Box_b heat heat equation and the Fourier transform
%
\subsection{The $\Box_b$-heat equation and the Fourier transform.} 
The heat equation defined on $(0,q)$-forms on $M$ is the initial value problem on $s\in(0,\infty)$ and
$(z,t)\in M$ given by 
\[
\begin{cases}
\displaystyle \frac{\p \rho}{\p s} +\Box_b\rho=0 &\text{in }(0,\infty)\times M \\
\rho(s=0,z,t) = \delta_0(z,t) &\text{on }\{s=0\}\times M
\end{cases}
\]
Here, $s$ is the time variable and $t\in \R^m$ is a spatial variable. Although we cannot find a closed form for  $\rho(s,z,t)$, we can find the partial Fourier
transform of $\rho(s,z,t)$ in the $t$-variables.

Given a variable
$\tilde t\in\R$, the \emph{(partial) Fourier transform} in $\tilde t$ is given by
\[
\hat f(\tau) = \frac{1}{\sqrt{2\pi}} \int_\R e^{-i\tilde t\tau} f(\tilde t)\, d\tilde t.
\]
If $f$ is a function of several variables $f(\tilde t_1, \dots, \tilde t_k)$ and, for example, we take the partial Fourier transform in $t_1$, we use the notation
$f ( \w{\tau}, \tilde t_2, \dots, \tilde t_k)$.

As we will see below, to compute the partial Fourier transform of $\rho(s,z,t) = \rho_s(z,t)$, it is enough to solve the (Fourier transform of the) $\Box_{LL}$-heat equation
\begin{equation}\label{eqn:Box-LL heat equation}
\begin{cases}
\displaystyle \frac{\p \rho}{\p s} +\Box_{LL}\rho=0 &\text{in }(0,\infty)\times M \\
\rho_{s=0}(z,t) = \delta_0(z,t) &\text{on }\{s=0\}\times M
\end{cases}
\end{equation}

We start by computing the partial Fourier transform in $t$ of $\Box_{LL}$,
denoted $\Box_{LL}^\lambda$.
We start with a reexamination of \eqref{eqn:BoxKK}. By taking the partial Fourier transform in $t$ of the formula for $\Box_{LL}$, the effect is to replace
$D_t$ with $i\lambda$. If we write $z = \sum_{k=1}^n z_k v_k$, then
\[
\Imm \Big\{ \sum_{k=1}^n \phi(z,v_k)\p_{z_k}\Big\} \cdot i\lambda = \sum_{k=1}^n \phi(v_k,v_k)\cdot i\lambda \Imm\Big\{ z_k\p_{z_k}\Big\}
= \sum_{k=1}^n i \mu^\lambda_k \Imm\Big\{ z_k\p_{z_k}\Big\}
\]
and
\[
 \sum_{k=1}^n \big( \phi(z,v_k) \cdot i\lambda \big) \big( \bar{\phi(z,v_k)} \cdot i\lambda \big)
 = \sum_{k=1}^n \big( z_k\phi(v_k,v_k) \cdot i\lambda \big) \big( \z_k\phi(v_k,v_k) \cdot i\lambda \big)
 = - \sum_{k=1}^n (\mu^\lambda_k)^2 |z_k|^2
\]
Consequently, \eqref{eqn:BoxKK} transforms to
\begin{equation}\label{eqn:Box-LL-lambda}
 \Box_{LL}^\lambda =  -\frac 14 \Delta
+ 2i \sum_{k=1}^n  \mu^\lambda_k \Imm\big\{ z_k\p_{z_k}\big\}
+ \sum_{k=1}^n  (\mu^\lambda_k)^2 |z_k|^2 \\
- \left( \sum_{k \in L} \mu^\lambda_k
- \sum_{k \not\in L}\mu^\lambda_k \right)
\end{equation}

We employ the following notation: for $1\leq j \leq {\nu(\lambda)}$, define $\eplj(L) = \eplj = \sgn(\mulj)$, if $j \in L$
and $\eplj= - \sgn(\mulj)$ if $j \not\in L$.

Our main result is the following.
%
%	Theorem: Computation of \rho(z,\w\lambda)
\begin{theorem}\label{thm:rho(z, lambda) computation}
For any $\lambda\in \R^m$, the partial Fourier transform to the fundamental solution to the $\Box_{LL}$-heat equation satisfies the heat equation
\[
\begin{cases}
\displaystyle \frac{\p \rho}{\p s} +\Box_{LL}^\lambda\rho=0 &\text{in }(0,\infty)\times \C^n \\
\rho(s=0,z,\w\lambda) = (2\pi)^{-m/2}\delta_0(z) &\text{on }\{s=0\}\times \C^n
\end{cases}
\]
and is given by
\[
\rho (s,x,y, \w \lambda) =
\frac{2^{n-{\nu(\lambda)}} (2 \pi)^{-(m/2+n)}} {s^{n-{\nu(\lambda)}}} e^{-\frac{|x''|^2+|y''|^2}s}  \prod_{j=1}^{\nu(\lambda)} \frac{2e^{s\eplj|\mulj|} \mulj} {\sinh (s\mulj)}
e^{-\mulj \coth(\mulj s) (x_j^2+y_j^2)}.
\]
\end{theorem}
Note that  $\mulj$ and $\coth (s\mulj)$ are real-valued and are odd in $\mulj$, so putting absolute values around the
$\mulj$ would not change the result. Therefore, there is Gaussian decay in $(x_j^2+y_j^2)$ for all $j$ when $\lambda\in\R^m$. Theorem \ref{thm:rho(z, lambda) computation}
generalizes Theorem 1.2 in \cite{BoRa09} in the case that $\tau\in\R$ and $\gamma = n-2q$.

We now cast the heat equation in terms of a weighted $\dbar$-problem in $\C^n$. Recall that $\bar Z_j = \frac{\p}{\p\z_j} + i\phi(z,v)\cdot D_t$. 
If we denote a superscript $\lambda$ for the partial Fourier transform in $t$, then
\[
\bar Z_j \mapsto \bar Z^\lambda_j = \frac{\p}{\p \z_j} -\phi(z,v_j)\cdot\lambda = e^{\phi(z,z)\cdot\lambda} \frac{\p}{\p\z_j} e^{-\phi(z,z)\cdot \lambda}.
\]
From the computation of $\bar Z_j$, the tangential Cauchy-Riemann operator $\dbar_b$ is defined on $(0,q)$-forms on $G$ by
\[
\dbar_b f(z) = \sum_{\atopp{K\in\I_{q+1}}{J\in\I_q}} \sum_{j=1}^n \eps^{jJ}_K \bar Z_j f_J(z)\, d\z^K
\]
where 
\[
\eps^{jJ}_K = \begin{cases} (-1)^\sigma &\text{if } \{j\}\cup J = K \text{ and } \sigma \text{ is the sign of the permutation taking } \{j\}\cup J \text{ to } K \\
0 &\text{otherwise}\end{cases}
\]
This means that if $g$ is a $(0,q)$-form in $\C^n$ and we treat $\lambda$ as a parameter, then the partial Fourier transform in $t$ of $\dbar_b$, denoted by 
$\dbar_b^\lambda$ is given by
\[
\dbar_b^\lambda g(z) = e^{\phi(z,z)\cdot\lambda} \dbar \{ e^{-\phi(z,z)\cdot\lambda} g\}
\]
where $\dbar$ is the usual Cauchy-Riemann operator on $\C^n$. Since $\Box_b = \dbar_b \dbar_b^* + \dbar_b^* \dbar_b$ where $\dbar_b^*$ is the $L^2$-adjoint
of $\dbar_b$, it follows that $\Box_b^\lambda = \dbar_b^\lambda \big(\dbar_b^\lambda\big)^* + \big(\dbar_b^\lambda\big)^* \dbar_b^\lambda$.
Thus, solving for the $\Box_b^{\lambda}$-heat kernel 
also yields the heat kernel associated to the weighted $\dbar$-problem on $\C^n$ with the weight $e^{-\phi(z,z)\cdot\lambda}$. 

% Corollary: weighted heat kernel on C^n
\begin{corollary} \label{cor:weighted heat kernel}
For any $\lambda\in\R^m$, the function
\[
H^\lambda(s,z,\tilde z) = (2\pi)^{m/2}\rho_s(z-\tilde z,\w\lambda) e^{-2i\lambda \cdot \Imm\phi(z,\tilde z)}
\]
satisfies the following: if
\[
H^\lambda\{f\}(s,z) = \int_{\C^n} H^\lambda(s,z,\tilde z) f(\tilde z)\, d\tilde z,
\]
then $H^\lambda\{f\}$ solves the initial value problem for the weighted heat equation:
\[
\begin{cases} 
(\p_s + \Box^\lambda_b)\{H^\lambda f\} =0 &\text{in } (0,\infty)\times\C^n \\
H^\lambda\{f\}(s=0, z) = f(z) &\text{on }\{s=0\}\times\C^n
\end{cases}
\]
In particular, the component of $H^\lambda(s,z,\tilde z)$ on $d\z^L$ for $L\in\I_q$ is
\[
H^\lambda_L(s,z,\tilde z) = \frac{2^{n-{\nu(\lambda)}} (2 \pi)^{-n}} {s^{n-{\nu(\lambda)}}} e^{-\frac{|z''-\tilde z''|^2}s} 
\prod_{j=1}^{\nu(\lambda)} \frac{2e^{s \eplj |\mulj|} \mulj } {\sinh (s \mulj )}
e^{- \mulj   \coth(\mulj s) |z_j-\tilde z_j|^2}
e^{-2i \lambda \cdot \Imm \phi(z,\tilde z)}.
\]
\end{corollary}
Note the formula for the heat kernel yields a standard Gaussian
solution for the Euclidean heat kernel in the zero eigenvalue
directions. Also, the disappearance of the $(2\pi)^{-m/2}$ owes to the fact that $\delta_0(z,\w\lambda) = (2\pi)^{-m/2}\delta_0(z)$.

%%%%%%%%%%%%%%%%%%%%%
%
%	SECTION:REPRESENTATION THEORY
%
%%%%%%%%%%%%%%%%%%%%%

\section{Representation Theory.}\label{sec:rep theory}

%
%	Irreducible unitary representations
%
\subsection{Irreducible unitary representations}
For $z=x+iy \in \C^n$, $t ,\lambda \in \R^m$,  and $\eta \in \C^{n-{\nu(\lambda)}}$, define
$\pi_{\lambda,\eta} (x,y,t) :L^2 (\R^{\nu(\lambda)}) \mapsto L^2(\R^{\nu(\lambda)})$
by
\[
\pi_{\lambda,\eta} (x,y,t) (h) (\xi)
= e^{i (\lambda \cdot t + 2\Rre(z''\cdot \bar\eta))} e^{-2i \sum_{j=1}^{\nu(\lambda)} \mulj y^\lambda_j ( \xi_j + x^\lambda_j)}
h (\xi +2x')
\]
for $h \in L^2 (\R^{\nu(\lambda)})$ (so $\xi \in \R^{\nu(\lambda)}$). Note that if $\eta = \zeta + i\varsigma$, then $\Rre(z''\cdot \bar\eta)) = x''\cdot \zeta + y''\cdot \varsigma$.

The map $\pi_{\lambda,\eta} (x,y,t)$ is unitary 
on $L^2 (\R^{\nu(\lambda)})$. Also, $\pi$ is a {\em representation }
for $G$, which means that for each $\lambda \in \Omega$,
$\pi_{\lambda,\eta}$ is a group homomorphism from $G$ to the group of 
unitary operators on $L^2(\R^{\nu(\lambda)})$. Verifying that $\pi_{\lambda,\eta}$ is a representation
is done in \cite{PeRi03}.

If $X$ is a right-invariant vector field, then $X$ gets
``transformed'' via $\pi_{\lambda,\eta}$ to an operator on $L^2(\R^{\nu(\lambda)})$ 
denoted by $T=d \pi_{\lambda,\eta} (X)$. This means
that
\begin{equation}
\label{rightinv}
X \{\pi_{\lambda,\eta} (g) \} =T \circ \pi_{\lambda,\eta} (g)
\end{equation}
as operators on $L^2(\R^{\nu(\lambda)})$. It is usually easy to identify
$T$ by seeing what happens at $g=0$ and using the right
invariance of $X$ to show that the above equation holds for all $g \in G$.
To clarify, let $R_g (g') = g' g$ and recall
that the vector field $X$
at the point $g$ is given by $X (g) = (R_g)_* \{X (0)\}$. If $
X \{\pi_{\lambda,\eta} \} ( 0 )  = T\circ \pi_{\lambda,\eta} (0)$,
then we have 
\begin{eqnarray*}
X \{\pi_{\lambda,\eta} \} (g)  &=&
(R_g)_* \{X (0) \} \{\pi_{\lambda,\eta} (g) \} \\
&=& X (g'=0) \{ \pi_{\lambda,\eta} (R_g(g'))\} \\
&=&X (g'=0) \{ \pi_{\lambda,\eta} (g') \pi_{\lambda,\eta} (g)\} 
\ \ \textrm{since} \ \pi \ \textrm{is a homomorphism} \\
&=&  \{ X (g'=0)\pi_{\lambda,\eta} (g' ) \} \circ \pi_{\lambda,\eta} (g) \\
&=& T\circ \pi_{\lambda,\eta} (g) 
\end{eqnarray*}
where the last equation uses the relationship of $X (g)$
and $\pi$ at $g=0$.

A similar computation shows that if $X^\ell$ is left invariant,
then 
\[
X^\ell \{\pi_{\lambda,\eta} \} (g)  =  \pi_{\lambda,\eta} (g) \circ T
\]
as operators on $L^2(\R^{\nu(\lambda)})$. Note that the order of $T$ and $\pi_{\lambda,\eta}$ is reversed
from (\ref{rightinv}). We will not dwell on this point as we prefer the use of right-invariant vector fields. 
The relationship $X \{\pi_{\lambda,\eta} \} (g)  = T \circ \pi_{\lambda,\eta} (g)$
is often expressed  using the shorthand:
$d \pi_{\lambda,\eta} (X) = T$.

From earlier,  we have the right invariant vector fields
\[
X_j=\partial_{v^\lambda_j} - 2 \Imm \phi (z,v^\lambda_j) \cdot D_t \\
\]
and
\[
Y_j =\partial_{Jv^\lambda_j} - 2 \Imm \phi (z,iv^\lambda_j) \cdot D_t = \partial_{Jv^\lambda_j} +2\Rre\phi (z,v^\lambda_j) \cdot D_t.
\]
where $J$ is the usual complex structure map on $\R^{2n} = \C^n$.
In view of (\ref{rightinv}),
we have the following relations: for 
\begin{align}
\label{righttransfereqn1}
 X_j \{\pi_{\lambda,\eta} \} (g) & = \begin{cases}  2 \p_{\xi_j} \circ \pi_{\lambda,\eta} (g)    & 1 \leq j \leq {\nu(\lambda)} \\
                               2i \zeta_j \circ \pi_{\lambda,\eta} (g) & {\nu(\lambda)}+1\leq j \leq n      \end{cases} \\
\label{righttransfereqn2}
Y_j \{\pi_{\lambda,\eta} \} (g)  &= \begin{cases} -2i \mulj \xi_j  \circ  \pi_{\lambda,\eta} (g)	& 1 \leq j \leq {\nu(\lambda)} \\
				2i \varsigma_j \circ \pi_{\lambda,\eta} (g) & {\nu(\lambda)}+1\leq j \leq n      \end{cases} \\
\label{righttransfereqn3}
\p_{t_k}\{\pi_{\lambda,\eta} \}(g)   &= i \lambda_k  \circ  \pi_{\lambda,\eta} (g) \qquad 1\leq k \leq m
\end{align}
as operators on $L^2(\R^{\nu(\lambda)})$. 
In the second equation, $\xi_j$ 
is thought of as a multiplication operator on $L^2(\R^{\nu(\lambda)})$,
i.e. $f(\xi) \mapsto f(\xi) \xi_j$. Equations
(\ref{righttransfereqn1}) and (\ref{righttransfereqn2})
are easily shown to hold 
at the origin since $X_j (0) = \p_{x_j}$
and $Y_j (0)= \p_{y_j}$, and the right invariance forces these
equations to hold at all $g \in G$.

Now we compute the ``transform'' of $\Box_{LK}$ in the coordinates $(z^\lambda_1, \dots, z^\lambda_n)$.
Note that 
\[
d \pi_{\lambda,\eta} [Z_j, \bar Z_\ell] =  \left\{
\begin{array}{cc}
-2 \mulj & \textrm{if} \ j = \ell \\
0 & \textrm{if} \ j \not= \ell.
\end{array}
\right.
\]
This follows from (\ref{righttransfereqn3}) and the fact
that the coordinates $(z^\lambda_1, \dots, z^\lambda_n)$
were chosen to diagonalize the form $\phi(z,\tilde z) \cdot \lambda$. In view of (\ref{boxLK}) and (\ref{righttransfereqn1}) -
(\ref{righttransfereqn3}), we have
\begin{equation}
d \pi_{\lambda,\eta} \Box_{LK} = \left\{
\begin{array}{cc}
- \Delta_\xi + |\eta|^2 +  \sum_{j=1}^{\nu(\lambda)} (\mulj)^2  \xi_j^2   -
\sum_{j=1}^{\nu(\lambda)} \eplj|\mulj| &
\textrm{if } \ K=L \\
0 & \textrm{if} \ K \not= L
\end{array}
\right.
\end{equation}
We will also need to transform the adjoint
of $\Box_{LK}$ which is defined as
\[
\int_{(z,t) \in G} \Box_{LK} \{ f(z,t) \} g(z,t) \, dx\, dy\, dt
=
\int_{(z,t) \in G} f(z,t) \Box^{\textrm{adj}}_{LK} \{g(z,t) \} \, dx\, dy\, dt
\]
(note: this is the ``integration by parts'' adjoint, not 
the $L^2$ adjoint, since there is no conjugation).
We have
\begin{equation}
Q^{\lambda, \eta, LK}_\xi := d \pi_{\lambda,\eta} \Box^{\textrm{adj}}_{LK}  = \left\{
\begin{array}{cc}
- \Delta_\xi + |\eta|^2 + \sum_{j=1}^{\nu(\lambda)} (\mulj)^2  \xi_j^2 +\sum_{j=1}^{\nu(\lambda)} \eplj|\mulj| &
\textrm{if } \ K=L \\
0 & \textrm{if} \ K \not= L
\end{array}
\right.
\end{equation}
(just a sign change for the last term on the right).
The subscript $\xi$ on $Q^{\lambda, \eta, LK}_\xi$ indicates 
that this is a differential operator in the $\xi$
variable (instead of the group variable $g=(x,y,t)$).
Below, we assume $L=K$ (otherwise the operator is zero)
and that $L$, $\lambda$, and $\eta$ are fixed. We drop the superscript $LL$ when its use is unambiguous.
In view of (\ref{boxLK}), and (\ref{righttransfereqn1})
through (\ref{righttransfereqn3}), we have
\begin{equation}
\label{keyeqn}
\Box^{\textrm{adj}}_{LL} \{ \pi_{\lambda,\eta} (g) \} = Q_\xi^{\lambda,\eta} \circ \pi_{\lambda,\eta} (g)
\end{equation}
as operators on $L^2(\R^{\nu(\lambda)})$.
We return to this key equation later.

%
%	Group Fourier Transform
%
\subsection{Group Fourier Transform.}
For $(z,t)\in G$, we express $(z,t)=(x,y,t) = (x',y',x'',y'',t) = (x',y', z'',t)$. The variable $z''$ may be thought of as in $\C^{n-{\nu(\lambda)}}$ or $\R^{2(n-{\nu(\lambda)})}$. 

For $f:G \mapsto \C$, we define the {\em group Fourier transform} of $f$
as the operator $T^{\lambda,\eta}_f: L^2(\R^{\nu(\lambda)}) \mapsto L^2(\R^{\nu(\lambda)})$ where for
$h \in L^2(\R^{\nu(\lambda)})$,
\begin{eqnarray*}
T^{\lambda,\eta}_f  \{h\} (\xi)
&=&  \int_{(z=x+iy,t) \in G} f(z,t) \pi_{\lambda,\eta} (z,t) (h)(\xi) \, dx\, dy\, dt \\
&=&  \int_{(z=x+iy,t) \in G} f(z,t)
e^{i (\lambda \cdot t + 2\Rre(z''\cdot\bar\eta))} e^{-2i \sum_{j=1}^{\nu(\lambda)} \mulj y^\lambda_j ( \xi_j + x^\lambda_j)}
h (\xi +2x')  \, dx\, dy\, dt.
\end{eqnarray*}
As before, $x_j$, $y_j$ are the coordinates
for $x,y \in \R^n$ relative to the basis $v_1^\lambda, \dots, v_n^\lambda$.
Note that 
\[
T^{\lambda,\eta}_f  \{h\}  (\xi) =  (2 \pi)^{(2n+m-{\nu(\lambda)})/2}
\int_{x' \in \R^{\nu(\lambda)}} f(x', 2 \w{\mu \circ (\xi +x')}, \w{-2\eta}, \w{-\lambda})h(\xi+2x') \, dx'.
\]
We have written $\mu^\lambda \circ (\xi +x')$
for $(\mu_1^\lambda (\xi_1 +x_1^\lambda), \dots, \mu_{\nu(\lambda)}^\lambda (\xi_{\nu(\lambda)} +x_{\nu(\lambda)}^\lambda))$.
We can also express $T^{\lambda,\eta}_f\{h\}$ as
\begin{multline}\label{keyfoureqn}
T^{\lambda,\eta}_f  \{h\}  (\xi)\\ 
=  (2 \pi)^{(2n+m-{\nu(\lambda)})/2}\int_{x' \in \R^{\nu(\lambda)}}\mathcal{F}_{x'',y,t} \{f (x,y,t) e^{-2i \sum_{j=1}^{{\nu(\lambda)}} \mulj
x_j y_j} \} (x', 2 \mu \circ \xi, -2\eta, - \lambda) h(\xi+2x') \, dx'.
\end{multline}
In the above notation, $\mathcal{F}_{x'',y,t}$ indicates the 
Fourier transform in the $(x'',y,t)$ variables only,
whereas $\mathcal{F}$ indicates the Fourier transform
in all variables (except $s$).

In view of (\ref{keyeqn}), we have
\begin{eqnarray}
Q^{\lambda,\eta}_\xi \{ T^{\lambda,\eta}_f  (h) (\xi) \}
&=& \int_{(z=x+iy,t) \in G} f(z,t) 
\Box^{\textrm{adj}}_{LL} \{ \pi_{\lambda,\eta} (z,t) h(\xi)  \}\, dx\, dy\, dt \nonumber \\
\label{kkeyeqn2}
&=& \int_{(z=x+iy,t) \in G} \Box_{LL} \{f(z,t) \}
\pi_{\lambda,\eta} (z,t) h(\xi) \, dx\, dy\, dt\,
\end{eqnarray} 

%%%%%%%%%%%%%%%%%%%%%
%
%	Heat Equations
%
%%%%%%%%%%%%%%%%%%%%%
\section{The Heat Equation.} \label{sec:heat equation and Group F.T.}
\subsection{The heat equation on $M$}
Our goal is to find a formula for the fundamental solution to the 
heat equation (\ref{eqn:Box-LL heat equation}).
%\begin{equation}
%\label{heateqn}
%\p_s \rho (s,z,t) = - \Box_{LL} \rho (s,z,t) \ \ 
%\textrm{with} \ \rho (s=0,z,t) = \delta_0 (z,t)
%\end{equation}
%where $\delta_0$ refers to the usual Dirac delta function
%supported at the origin. 
We know abstractly that $\rho$ exists: $\Box_{LL}$ is self-adjoint and nonnegative, so $e^{-s\Box_{LL}}$ is a well-defined, bounded linear operator
on $L^2(G)$ with norm at most 1. It has an integral kernel by the Riesz Representation Theorem, and the existence of $\rho$ follows easily. See
\cite{Rai06h} for details.

Let us apply the group Fourier transform to $\rho$ and recall that $\rho_s(z,t) = \rho(s,z,t)$.
Define the operator $U^{\lambda,\eta} (s): L^2(\R^{\nu(\lambda)}) \mapsto L^2(\R^{\nu(\lambda)})$ by
\begin{equation}
\label{Ulambda}
U^{\lambda,\eta} (s) \{h\}  (\xi) = T_{\rho_s}^{\lambda,\eta} \{h\}(\xi)
=  \int_{(z,t) \in G} \rho_s (z,t) \pi_{\lambda,\eta} (z,t) h(\xi) \, dz\, dt.
\end{equation}
In view of (\ref{kkeyeqn2}) and the fact that 
$\rho_s (z,t)$ solves the heat equation, we have
\begin{eqnarray*}
Q_\xi^{\lambda,\eta} \{ U^{\lambda,\eta}(s) \{h\} (\xi) \}
&=& \int_{(z,t) \in G} \Box_{LL} \{ \rho_s (z,t) \}
\pi_{\lambda,\eta} (z,t) h(\xi) \, dz\, dt \\
&=& - \p_s \left\{ \int_{(z,t) \in G}  \rho_s (z,t) 
\pi_{\lambda,\eta} (z,t) h(\xi) \, dz\, dt \right\} \\
&=& -\p_s \big\{ U^{\lambda,\eta}(s) \{h\} (\xi) \big\}
\end{eqnarray*}
Also
\[
U^{\lambda,\eta}(s=0)\{h\}(\xi) =  T_{\delta_0 }^{\lambda,\eta} \{h\}(\xi) = h(\xi).
\]
Therefore, we conclude that $U^{\lambda,\eta}(s) $ satisfies
the following boundary value problem:
\begin{equation}
\label{hermitediffeq}
Q_\xi^{\lambda,\eta} \{ U^\lambda (s) \} = -\p_s \{ U^\lambda (s) \} 
\ \ \textrm{and} \ \ U^\lambda (s=0) = \textrm{Id}
\end{equation}
where $\textrm{Id}$ is the identity operator on $L^2(\R^{\nu(\lambda)})$.
This is a Hermite equation similar to, though more complicated than, the one we solved
in the Heisenberg group case \cite{BoRa09}.
So, our approach is to proceed as follows: \textbf{1)} 
explicitly solve this Hermite equation, and then \textbf{2)}  recover
the fundamental solution to the heat equation.

As to the second task, we let $a \in \R^{\nu(\lambda)}$ be an arbitrary
vector, and then define $h_a (\xi)= (2\pi)^{-n-m/2} e^{- i\xi \cdot a}$. Let
\[
u^{\lambda,\eta}(s,a, \xi)  = U^{\lambda,\eta}(s) \{h_a\}(\xi) .
\]
The above definition needs explanation since $h_a \not\in L^2(\R^{\nu(\lambda)})$. For each fixed $s>0$,  $\rho_s\in L^2(G)$ and we can approximate $\rho_s$ by 
$\rho_s^\delta\in L^1\cap L^2(G)$ (e.g., by multiplying $\rho_s$ with an appropriate test function). Then, as we see below, we can define
$U^{\lambda,\eta}_\delta(s)\{ h_a\}(\xi) =  T_{\rho_s^\delta}^{\lambda,\eta} \{h_a\}(\xi)$ since 
in view of equation (\ref{keyfoureqn}),
\begin{align*}
&U^{\lambda,\eta}_\delta (s) \{h_a\}(\xi)\\
&= \frac{1}{(2 \pi)^{{\nu(\lambda)}/2} }\int_{x'\in\R^{{\nu(\lambda)}}} \mathcal{F}_{x'',y,t} \left\{ \rho_s^\delta(x,y,t)
e^{-2i \sum_{j=1}^{{\nu(\lambda)}} \mulj y_j x_j} \right\} (x', 2 \mu \circ \xi, -2\eta, -\lambda) e^{-i ( \xi + 2x') \cdot a} \, dx' \\
&=\mathcal{F} \left\{ \rho_s^\delta (x,y,t)
e^{-2i \sum_{j=1}^{\nu(\lambda)} \mulj y_j x_j} \right\}
(2a, 2 \mu^\lambda \circ \xi, -2\eta, -\lambda) e^{-i \xi \cdot a}
\end{align*}
By the definition of the Fourier transform in $L^2$, 
\begin{multline*}
\mathcal{F} \left\{ \rho_s^\delta (x,y,t)
e^{-2i \sum_{j=1}^{\nu(\lambda)} \mulj y_j x_j} \right\}
(2a, 2 \mu^\lambda \circ \xi, -2\eta, -\lambda) e^{-i \xi \cdot a} \\
\longrightarrow 
\mathcal{F} \left\{ \rho_s (x,y,t)
e^{-2i \sum_{j=1}^{\nu(\lambda)} \mulj y_j x_j} \right\}
(2a, 2 \mu^\lambda \circ \xi, -2\eta, -\lambda) e^{-i \xi \cdot a}
\end{multline*}
in $L^2(\R^{\nu(\lambda)})$ as $\delta\to 0$. Thus, $u^{\lambda,\eta}(s,a,\xi)$ is well-defined.
In the above computation, we view $\eta = (\zeta,\varsigma)\in \R^{2(n-{\nu(\lambda)})}$. Also,
the motivation for the choice of $h = h_a$ is that 
it offers the ``missing'' exponential needed to 
relate the full Fourier transform of $\rho_s$ with $u^{\lambda,\eta}$.
Now it is just a matter of unraveling the equation
\begin{equation}
\label{ulambda1}
u^{\lambda,\eta}(s,a, \xi) = 
\mathcal{F} \left\{ \rho_s (x,y,t)
e^{-2i \sum_{j=1}^{\nu(\lambda)} \mulj y^\lambda_j x^\lambda_j} \right\}
(2a, 2 \mu^\lambda \circ \xi, -2\eta, -\lambda) e^{-i \xi \cdot a}
\end{equation}
for $\rho_s$ using the inverse Fourier transform.

Before we go on, let us remark that had we used
left invariant vector fields rather than right invariant
ones, then the transformed operator, $Q_\xi^{\lambda,\eta}$
would appear on the right of the group transform.
That is to say, we would be trying to solve 
the following analogue of (\ref{hermitediffeq})
\[
 \p_s \{T^\lambda_{\rho_s}  \}
= -T^\lambda_{\rho_s} \tilde Q_\xi^{\lambda,\eta} \ \ 
\textrm{and} \ \ T^\lambda_{\rho_{s=0}} = \textrm{Id}
\]
where $\tilde Q_\xi^{\lambda,\eta}$ is a Hermite type differential
operator similar to $ Q_\xi^{\lambda,\eta}$.
Note the transform operator $T^\lambda$  is now intertwined
with the differential operators (i.e. $\p_s$ is on the 
left side and $\tilde Q_\xi^{\lambda,\eta}$ is on the right).
Since the inversion formula for the group transform operator
is complicated (see  \cite{PeRi03}), it would appear
that this method of using left invariant vector fields
is more difficult to unravel a formula for $\rho$.
%But I may be overlooking something simple.

\subsection{Weighted heat equation.} 
Our objective is to compute $\rho_s(x',y',\w\eta,\w\lambda)$ by
solving the weighted heat equation obtained by taking the partial Fourier transform in  the $t$ and $(x'',y'')$-variables.
We obtain the $\Box^{\lambda,\eta}_{LL}$-heat equation
\[
\begin{cases} \p_s \rho_s(x',y', \w\eta, \w{\lambda}) = -\Box^{\lambda,\eta}_{LL} \rho_s(x',y', \w\eta, \w{\lambda}) \\
\rho_{s=0}(x',y', \w\eta, \w{\lambda}) =(2 \pi)^{-m/2-(n-{\nu(\lambda)})}\delta_0 (x',y')
\end{cases}
\]
From (\ref{eqn:Box-LL-lambda}), we have
\[
\Box^{\lambda,\eta}_{LL} = -\frac 14\Delta  + \frac 14 |\eta|^2 + 2i \sum_{j=1}^{\nu(\lambda)}  \mulj
\Imm \{ z_j  \p_{z_j} \} +  \sum_{j=1}^{\nu(\lambda)} |z_j \mulj|^2 - \sum_{j=1}^{\nu(\lambda)} \epj |\mulj|.
\] 
In the following computation, we find a formula of $\rho_s(x',y',\w\eta,\w\lambda)$. 
Observe that
$\mu^{-\lambda}_j = -\mulj$ and $\eps^{-\lambda}_j = -\eplj$. (Note that $v^{-\lambda}_j = v^{\lambda}_j$, so we can continue to suppress the $\lambda$ superscript
on $x_j$ and $y_j$.)
We unravel  (\ref{ulambda1}) to obtain (with $a,b\in \R^{{\nu(\lambda)}}$)
\begin{equation}
\label{rhos}
\rho_s (x',y', \w\eta, \w\lambda) = 
e^{-2i \sum_{j=1}^{{\nu(\lambda)}}  \mulj x_j y_j}
\mathcal{F}^{-1}_{a,b} \left(e^{-\frac i4 \sum_{j=1}^{{\nu(\lambda)}} a_j b_j /\mulj}
\tilde u^{\lambda,\eta} (s, a, b) \right) (x',y')
\end{equation}
where  $\tilde u^{\lambda,\eta} (s, a, b) = u^{-\lambda, -\frac 12\eta} (s, a/2, b/(2 \mu^{- \lambda}))$
and $b/(2 \mu^{-\lambda})$ is the vector quantity whose
$j$th component is $b_j/(2 \mu_j^{-\lambda})$.
As we shall see, the inverse Fourier transform in the $a$, and $b$ variables
will be relatively simple (using Gaussian integrals).
In the next section, we use Hermite functions to solve for $\tilde u^{\lambda,\eta}$
on the ``transform'' side. Then we return to the above formula to compute
$\rho_s (x',y', \w\eta, \w \lambda)$.

%%%%%%%%%%%%%%%%%%%%%%%%%
%
%	SOLVING THE HEAT EQUATIONS	
%
%%%%%%%%%%%%%%%%%%%%%%%%%
\section{Computing the heat kernels}\label{sec:heat kernel computation}
In this section, we prove Theorem \ref{thm:rho(z, lambda) computation} and Corollary \ref{cor:weighted heat kernel}.

\subsection{Hermite Functions.} Our starting point is 
equation (\ref{hermitediffeq}), which we restate as
$Q_\xi^{\lambda,\eta} \{ U^{\lambda,\eta}(s) \} = -\p_s \{ U^{\lambda,\eta}(s) \}$
where 
\[
Q_\xi^{\lambda,\eta} = - \Delta_\xi  + |\eta|^2 + \sum_{j=1}^{\nu(\lambda)}  (\mulj \xi_j)^2  +
\sum_{j=1}^{\nu(\lambda)} \epj | \mulj|.
\]
We use Hermite functions to solve this equation. 
For a nonnegative integer $\ell$, define
\[
\psi_\ell(x) = 
\frac{(-1)^\ell} {2^{\ell/2} \pi^{1/4}(\ell!)^{1/2}}\frac{d^\ell}{dx^\ell}\{ e^{-x^2} \} e^{x^2/2} ,
\ \ x \in \R.
\]
Each $\psi_\ell$ has unit $L^2$-norm on the real line and satisfies
the equation
\[
- \psi_\ell ''(x) + x^2 \psi_\ell (x) = (2\ell+1) \psi_\ell(x),
\]
see \cite{Tha93}, (1.1.9).
For $\lambda\in\R^m\setminus\{0\}$, define
\[
\psi^\lambda_{\ell_j} (\xi_j) = \psi_{\ell_j} (|\mulj|^{1/2} \xi_j)
|\mulj|^{1/4}.
\]
Each $\psi^\lambda_{\ell_j} (\xi_j)$ has unit $L^2$-norm on $\R$
and hence $\psi^\lambda_\ell$ has unit $L^2$-norm on $\R^{\nu(\lambda)}$.
An easy calculation shows that 
\begin{equation}
\label{hermite1}
(-\p_{\xi_j \xi_j} + (\mulj\xi_j)^2) \{ \psi_{\ell_j}^\lambda (\xi_j) \}
= (2\ell_j+1) \psi_{\ell_j}^\lambda (\xi_j) |\mulj|.
\end{equation}

For $s>0$, we claim that $U^{\lambda,\eta}(s) : L^2 (\R^{\nu(\lambda)}) \mapsto L^2(\R^{\nu(\lambda)}) $ as defined in
(\ref{Ulambda}) is given by
\[
U^{\lambda,\eta}(s) = e^{-s|\eta|^2} \bigotimes_{j=1}^{\nu(\lambda)}
\sum_{\ell_j=0}^\infty e^{-[(2\ell_j+1) + \epj]|\mulj| s }
P^\lambda_{\ell_j}
\]
where $P^\lambda_{\ell_j} $ is the $L^2$ projection of a smooth
function of polynomial growth in the variable $\xi_j$ onto the 
space spanned by $\psi^\lambda_{\ell_j} (\xi_j)$, and where 
$\bigotimes_{j=1}^{\nu(\lambda)}$ is the tensor product (so that the output
of $U^{\lambda,\eta}(s)$ is a function of $\xi_1, \dots , \xi_{\nu(\lambda)}$).
For shorthand, we write
\[
E^\lambda_{\ell_j} (s)=  e^{-[(2\ell_j+1) + \epj]|\mulj| s }.
\]
We then have $U^{\lambda,\eta}(s) = e^{-s|\eta|^2} \bigotimes_{j=1}^{\nu(\lambda)} \sum_{\ell_j=0}^\infty E^\lambda_{\ell_j} (s) P^\lambda_{\ell_j} $.
Using the product rule, we compute
\begin{align*}
&\p_s \{U^{\lambda,\eta}(s) \} = - |\eta|^2 U^{\lambda,\eta}(s) + e^{-s|\eta|^2}\sum_{j=1}^{\nu(\lambda)} \p_s \left( \sum_{\ell_j=0}^\infty E^\lambda_{\ell_j} (s) P^\lambda_{\ell_j} \right)
\bigotimes_{\atopp{k=1}{k \neq j}}^{\nu(\lambda)} \sum_{\ell_k=0}^\infty E^\lambda_{\ell_k} (s) P^\lambda_{\ell_k} \\
&=- |\eta|^2 U^{\lambda,\eta}(s) + e^{-s|\eta|^2}\sum_{j=1}^{\nu(\lambda)}\sum_{\ell_j =0}^\infty -[(2\ell_j+1)+\epj]|\mulj|   
e^{-[(2\ell_j+1) + \epj]|\mulj| s } P^\lambda_{\ell_j}
\bigotimes_{\atopp{k=1}{k \neq j}}^{\nu(\lambda)} \sum_{\ell_k=0}^\infty E^\lambda_{\ell_k} (s) P^\lambda_{\ell_k} \\
&=- |\eta|^2 U^{\lambda,\eta}(s) + e^{-s|\eta|^2}\sum_{j=1}^{\nu(\lambda)} \sum_{\ell_j=0}^\infty  \left(\p_{\xi_j \xi_j} - (\mulj\xi_j)^2 -
\epj|\mulj|\right)  \circ 
e^{-[(2\ell_j+1) + \epj]|\mulj| s } P^\lambda_{\ell_j}
\bigotimes_{\atopp{k=1}{k \neq j}}^{\nu(\lambda)} \sum_{\ell_k=0}^\infty E^\lambda_{\ell_k} (s) P^\lambda_{\ell_k}
\end{align*}
where the last equality uses (\ref{hermite1}).
Since the differential operator on the right is independent of $\ell_j$,
we can factor it to the left of $\sum_{\ell_j}$ to obtain
\[
\p_s \{ U^{\lambda,\eta}(s) \} = - Q_\xi^{\lambda,\eta} \{U^{\lambda,\eta}(s) \}.
\]
Since the Hermite functions, $\psi^\lambda_\ell$, form an orthonormal 
basis for $L^2(\R)$, $U^{\lambda,\eta}(s=0)$ is just the identity
operator. Thus $U^{\lambda,\eta}(s)$ solves (\ref{hermitediffeq}).

As above, we
apply $U^{\lambda,\eta}(s)$ to the function 
$h_a (\xi)= (2 \pi)^{-n-m/2}e^{-i \xi \cdot a}$ to obtain the fundamental solution $\rho_s$.
We therefore obtain
\begin{align*}
u^{\lambda,\eta}(s, a, \xi) &= U^{\lambda,\eta}(s) \{ h_a (\xi) \} \\
&= (2 \pi)^{-n-m/2} e^{-s|\eta|^2} \prod_{j=1}^{\nu(\lambda)} \sum_{\ell_j=0}^\infty E^\lambda_{\ell_j}(s) P^\lambda_{\ell_j} \{e^{- i \xi_j  a_j} \}
\end{align*}
Since $h_a$ belongs to $L^\infty (\R^{\nu(\lambda)})$ and {\em not} in $L^2 (\R^{\nu(\lambda)})$, the above sum converges
\emph{a priori} in the sense of tempered distributions (as opposed to $L^2$ convergence). Earlier, we argued that we can obtain $U^{\lambda,\eta}(s)\{h_a\}$ via a standard
approximation argument, 
however, we will see below that the convergence is much stronger and the result
is a smooth function in $s, \ a, \ \xi$. 
Each projection term on the right is
\begin{align*}
 P^\lambda_{\ell_j} (e^{- i  \xi_j  a_j} ) 
 &= \left(  \int_{\tilde \xi_j \in \R} e^{-i \tilde \xi_j a_j}  \psi_{\ell_j} ( |\mulj|^{1/2} \tilde \xi_j)  |\mulj|^{1/4} \, d \tilde \xi_j \right) |\mulj|^{1/4} 
 \psi_{\ell_j}( |\mulj|^{1/2} \xi_j) \\
 &= (2 \pi)^{1/2}\w{\psi_{\ell_j}} (a_j/|\mulj|^{1/2}) \psi_{\ell_j} ( |\mulj |^{1/2} \xi_j) \\
 &= (2 \pi)^{1/2}(-i)^{\ell_j} 
 \psi_{\ell_j} (a_j/|\mulj|^{1/2}) \psi_{\ell_j} ( |\mulj|^{1/2} \xi_j)
\end{align*}
where the last equality uses a standard fact about Hermite functions
that they equal their Fourier transforms up to a factor of 
$(-i)^{\ell_j}$. Substituting this expression on the right into the definition
of $u^{\lambda,\eta}(s, a, \xi)$, we obtain
\[
u^{\lambda,\eta} (s, a, \xi) =
(2 \pi)^{-n-m/2 +{\nu(\lambda)}/2} e^{-s|\eta|^2} \prod_{j=1}^{\nu(\lambda)} \sum_{\ell_j=0}^\infty E^\lambda_{\ell_j} (s) (-i)^{\ell_j}
\psi_{\ell_j} (a_j/|\mulj|^{1/2}) \psi_{\ell_j} ( |\mulj|^{1/2} \xi_j).
\]
This function satisfies
\begin{align*}
\p_s u^{\lambda,\eta} (s, a, \xi) &= - Q^{\lambda,\eta}_\xi \{u^{\lambda,\eta} (s, a , \xi) \} \\
u^{\lambda,\eta} (s=0, a, \xi) &= h_a (\xi) = (2 \pi )^{-n-m/2} e^{-i a \cdot \xi}.
\end{align*}

In view of (\ref{rhos}), for computing $\rho_s (x',y', \w\eta, \w \lambda)$,
we need to compute
\[
\tilde u^{\lambda,\eta} (s, a, b) = u^{-\lambda,-\frac 12\eta} (s, a/2, b/(2 \mu^{- \lambda}))
\]
where $b/(2 \mu^{-\lambda})$ is the vector quantity whose
$j$th component is $b_j/(2 \mu^{-\lambda}_j)$.
From the previous equality, and using that $\mu^{-\lambda}_j = - \mulj$,
$\eps^{-\lambda}_j = - \eplj$,
we have
\begin{multline*}
\tilde u^{\lambda,\eta} (s, a, b) =
(2 \pi)^{-\frac 12(n+m + (n-{\nu(\lambda)}))} e^{-s \frac{|\eta|^2}4} \prod_{j=1}^{\nu(\lambda)} e^{-(1 - \epj) | \mulj| s} \\
\times \sum_{\ell_j=0}^\infty (-i)^{\ell_j}
\psi_{\ell_j} (a_j/2|\mulj|^{1/2}) \psi_{\ell_j} (b_j |\mulj|^{1/2}/2 \mu^{-\lambda}_j)
e^{-2\ell_j |\mulj|s}.
\end{multline*}
Let
\begin{equation}
\label{shorthand}
S_j = e^{-2 |\mulj|s}, \ \alpha_j = \frac{a_j} {2 |\mulj|^{1/2}}, \ \ 
\beta_j = \frac {-b_j |\mulj|^{1/2}} {2 \mulj}.
\end{equation}
Then
\[
\tilde u^{\lambda,\eta} (s, a, b) = (2 \pi)^{-\frac 12(n+m + (n-{\nu(\lambda)}))}e^{-s\frac {|\eta|^2}4}
\prod_{j=1}^{\nu(\lambda)}  S_j^{(1-\epj)/2}
\sum_{\ell=0}^\infty (-iS_j)^\ell \psi_\ell( \alpha_j) \psi_\ell ( \beta_j).
\]
Using Mehler's formula (\cite{Tha93}, Lemma 1.1.1) for Hermite functions, we obtain
\[
\tilde u^{\lambda,\eta} (s, a, b) =
(2 \pi)^{-(m/2+n)} 2^{{\nu(\lambda)}/2} e^{-s\frac {|\eta|^2}4} \prod_{j=1}^{\nu(\lambda)}  S_j^{(1-\epj)/2} \frac 1{\sqrt{1+S_j^2}}
e^{- \frac 12 \big(\frac{1-S_j^2}{1+S_j^2}\big) (\alpha_j^2 + \beta_j^2)
- \frac{2i S_j \alpha_j \beta_j}{1+S_j^2}}.
\]
The series for $\tilde u^{\lambda,\eta}$ converges 
in $C^\infty$ on the unit disk in $\C$, and therefore the series for $\tilde u^{\lambda,\eta}$ converges in $\C^\infty$ for $s>0$, justifying
many previous computations (which held \emph{a priori} in the 
category of distributions).
%
%	Finish the proof of the weighted heat kernel computation
%
\subsection{Finishing the proof of Theorem \ref{thm:rho(z, lambda) computation}.}
In view of Equation (\ref{rhos}), to determine $\rho_s(x',y', \w\eta, \w \lambda)$,
we must compute
\[
\mathcal{F}^{-1}_{a,b} \left(e^{-i \sum_{j=1}^{\nu(\lambda)} a_j b_j /(4 \mulj)}
\tilde u^{\lambda,\eta} (s, a, b) \right) (x',y').
\]
Using (\ref{shorthand}), and simplifying we obtain 
\[
e^{-i \sum_{j=1}^{\nu(\lambda)} a_j b_j /(4 \mulj)}
\tilde u^{\lambda,\eta} (s, a, b)
=(2 \pi)^{-(m/2+n)}e^{-s\frac {|\eta|^2}4} \prod_{j=1}^{\nu(\lambda)} \frac {e^{\epj|\mulj| s} } {\sqrt{ \cosh(2 |\mulj|s)}}
e^{-A_j (a_j^2+b_j^2)/2 - i B_j a_j b_j}
\]
where
\[
A_j = \frac{\tanh (2 |\mulj| s)} { 4 |\mulj|}, \ \ 
B_j = \frac { \sinh^2(|\mulj|s)} { 2 \mulj \cosh (2 |\mulj| s)}.
\]
After an exercise in computing Gaussian integrals, we obtain
\begin{multline*}
\mathcal{F}^{-1}_{a,b} \left\{ e^{-i \sum_{j=1}^{\nu(\lambda)} a_j b_j /(4 \mulj)}
\tilde u^{\lambda,\eta} (s, a, b) \right \} (x',y')\\ 
= (2 \pi)^{-(m/2+n)}e^{-s\frac {|\eta|^2}4} \prod_{j=1}^{\nu(\lambda)}   \frac{e^{\epj|\mulj| s}} {\sqrt{ \cosh(2 |\mulj|s)}}
 \frac{e^{\frac {-A_j} {2(A_j^2+B_j^2)}(x_j^2+y_j^2) + i \frac{B_j x_j y_j} {A_j^2+B_j^2}}} {\sqrt{A_j^2+B_j^2}}.
\end{multline*}
After simplifying,
\[
\frac {-A_j} {2(A_j^2+B_j^2)} =-\mulj A_j /B_j, 
\ \ \ \frac {B_j} {A_j^2+B_j^2} = 2 \mulj, \ \ 
\sqrt{\cosh(2|\mulj| s)} \sqrt{A_j^2+B_j^2}
=\frac{\sinh(s|\mulj|)} {2|\mulj|}.
\]
The previous expression becomes
\begin{multline*}
\mathcal{F}^{-1}_{a,b} \left(e^{-i \sum_{j=1}^{\nu(\lambda)} a_j b_j /(4 \mulj)}
\tilde u^{\lambda,\eta} (s, a, b) \right) (x',y')  \\
= (2 \pi)^{-(m/2+n)} e^{-s\frac {|\eta|^2}4} \prod_{j=1}^{\nu(\lambda)} 
\frac{2e^{\epj|\mulj| s }|\mulj|} {\sinh (s|\mulj|)}
e^{-\mulj (A_j/B_j) (x_j^2+y_j^2) +2i\mulj x_j y_j}.
\end{multline*}
In view of (\ref{rhos}), the fundamental solution $\rho_s (x',y', \w\eta, \w \lambda)$
to the weighted heat equation is obtained by multiplying this expression
by $\prod_{j=1}^{\nu(\lambda)} e^{-2i \mulj x_j y_j}$ which cancels
the similar expression on the right side. We therefore obtain
\[
\rho_s (x',y', \w\eta, \w \lambda) =
(2 \pi)^{-(m/2+n)}e^{-s\frac {|\eta|^2}4}  \prod_{j=1}^{\nu(\lambda)} \frac{2e^{\epj|\mulj| s} |\mulj|} {\sinh (s|\mulj|)}
e^{-\mulj (A_j/B_j) (x_j^2+y_j^2)}.
\]
Note that the rightmost exponent can be rewritten as
\[
-\mulj (A_j/B_j) (x_j^2+y_j^2)=
- \frac{|\mulj|  \sinh(2 |\mulj| s)} {2 \sinh^2 (|\mulj|s)} (x_j^2+y_j^2) = -\mulj \coth(\mulj s)(x_j^2+y_j^2).
\]
Consequently,
\[
\rho_s (x,y, \w \lambda) =
\frac{2^{n-{\nu(\lambda)}} (2 \pi)^{-(m/2+n)}} {s^{n-{\nu(\lambda)}}} e^{-\frac{|x''|^2+|y''|^2}s}  \prod_{j=1}^{\nu(\lambda)} \frac{2e^{\epj|\mulj| s} |\mulj|} {\sinh (s|\mulj|)}
e^{-\mulj \coth(\mulj s)(x_j^2+y_j^2)}.
\]
This completes the proof of Theorem \ref{thm:rho(z, lambda) computation} for $\lambda\in\Omega$.

%
% Proof of the corollary
%

\subsection{The proof of Corollary \ref{cor:weighted heat kernel}}
In this subsection and the next, we show that the following kernel:
\begin{align}
\label{wtheatker}
&H^\lambda (s, z,\tilde z) = (2\pi)^{m/2}\rho_s (z-\tilde z, \w \lambda)
e^{-2i \lambda \cdot \Imm \phi(z,\tilde z)}\\
&= \frac{2^{n-{\nu(\lambda)}} (2 \pi)^{-n}} {s^{n-{\nu(\lambda)}}} e^{-\frac{|z''-\tilde z''|^2}s} 
\prod_{j=1}^{\nu(\lambda)} \frac{2e^{\epj|\mulj| s} |\mulj|} {\sinh (s|\mulj|)}
e^{-\mulj \coth(\mulj s) |z_j-\tilde z_j|^2}
e^{-2i \lambda \cdot \Imm \phi(z,\tilde z)} \nonumber
\end{align}
is the heat kernel for the weighted $\Box_b$. Here, $z=x+iy$ and $\tilde z=\tilde x+i\tilde y$.
Note that $H$ is conjugate
symmetric, i.e. $H^\lambda(s, \tilde z, z) = \bar{H^\lambda (s, z, \tilde z)}$.
We will show that the heat kernel has the following properties: if 
$f \in L^2(\C^n)$, then
\[
H^\lambda \{f\}(s,x,y) = \int_{\R^n \times \R^n}
H^\lambda (s, x,y, \tilde x, \tilde y) f(\tilde x,\tilde y) \, d\tilde x\, d\tilde y
\]
is the solution to the following boundary value problem for the heat equation:
\[
(\p_s +\Box^\lambda_b) \{H^\lambda f\} = 0, \ \ \ 
H^\lambda \{f\}(s=0, x,y) = f(x,y).
\]

%
%	Group Convolution and Twisted Convolution
%
\subsection{Group Convolution and Twisted Convolution}
To motivate the above formula,  we consider the fundamental solution
to the (full) unweighted heat equation: $\rho_s (x,y,t)$. For 
a function $f_0 \in L^2 (\C^n \times \R^m)$, and $g=(z,t) \in \C^n \times \R^m$ define
\begin{equation}
\label{heatker}
H\{f_0\} (s,g)= (\rho_s* f_0)(g) = \int_{\tilde g}  \rho_s (g[\tilde g]^{-1}) f_0(\tilde g) \, d\tilde g
\end{equation}
where $*$ is the group convolution and $g[\tilde g]^{-1}$ is the group multiplication
of $g$ by the inverse of $\tilde g$. If $X$ is a right invariant vector field,
then 
\[
XH\{f_0\}(s, g)  = \int_{\tilde g}  (X\rho_s) (g[\tilde g]^{-1}) f_0(\tilde g) \, d\tilde g.
\]
Since $\Box_b$ is comprised of right invariant vector fields and 
$\rho_s$ satisfies the $\Box_b$-heat equation, we therefore have
\[
(\p_s +\Box_b)\{H(f_0)\} =0.
\]
In addition, the following initial condition holds:
\[
H\{f_0\}(s=0,g) = \int_{\tilde g}  \rho_{s=0} (g[\tilde g]^{-1})f_0(\tilde g) \, d\tilde g =
f_0(g) 
\]
since $\rho_{s=0} (z,t)$ is the Dirac delta function centered at $(z,t)=0$.

Note that $H^\lambda \{f\}(s,x,y) = (2\pi)^{m/2}H\{ f\}(s,x,y, \w \lambda)$,
which is the partial Fourier transform in the $t$ variable of 
$H\{f\}(s,x,y,t)$.
We will now show the
Fourier transform in the $t$-variable transforms the 
group convolution to a ``twisted convolution'', which we now define.
Suppose $F$ and $G$ are in $L^2(\C^n)$, and $\lambda \in \R^m$.
Following Stein \cite{Ste93}, page 552, we let
\[
(F*_\lambda G)(z) = \int_{\tilde z \in \C^n} F(z-\tilde z)G(\tilde z) 
e^{-2i \lambda \cdot \Imm \phi (z,\tilde z)} \, d\tilde z.
\]
The arguments in \cite{Ste93}, page 552 with $\langle z, \tilde z \rangle$
replaced by $2 \Imm \phi (z,\tilde z)$ show the following: if $F_0, G_0 \in L^2(\C^n\times\R^m)$, then
\[
(F_0 * G_0)(z, \w \lambda) = (2\pi)^{m/2}
(F_0( \cdot, \w \lambda ) *_\lambda G_0( \cdot, \w \lambda)) (z).
\]
Now suppose $f \in L^2(\C^n)$ is given
and let $f_0 (z, t)= (2 \pi)^{m/2} f(z) \delta_0(t)$,
so that $f_0 (z, \w \lambda)  = f(z)$.
With $H^\lambda$ given as in (\ref{wtheatker}), we can take the partial Fourier transform in $t$ 
of (\ref{heatker}) and use the above relationship
to obtain
\begin{eqnarray*}
H^\lambda (f)(s, z) & =& \int_{\tilde z \in \C^n}  H^\lambda (s, z, \tilde z) f(\tilde z) \, d\tilde z \\
&=& \int_{\tilde z \in \R^n, \tilde y \in \R^n} (2\pi)^{m/2}
\rho_s( x-\tilde x, y-\tilde y, \w \lambda) f(\tilde x,\tilde y) 
e^{-2i \lambda \cdot \Imm \phi (z,\tilde z)} \, d\tilde x\, d\tilde y \\
&=& (2\pi)^{m/2}(\rho_s (\cdot, \w \lambda) *_\lambda f_0 (\cdot , \w \lambda) )(z)\\
&=& (\rho_s * f_0)(z,\w \lambda) \\
&=& H(f_0)(s,z, \w \lambda)
\end{eqnarray*}
Since $H(f_0)$ satisfies the $\Box_b$-heat equation,
$H (f_0)(s,z, \w \lambda) = H^\lambda (f)(s,z)$ satisfies
the weighted heat equation, i.e. 
\[
(\p_s + \Box^\lambda_b ) \{ H^\lambda (f) \} =0.
\]
The initial condition $H^\lambda (f) (s=0, z) = f(z)$ is also 
satisfied because
\begin{eqnarray*}
H^\lambda (f) (s=0, z) &=& H(f_0)(s=0,z, \w \lambda) \\
&=& f_0(z, \w \lambda) \\
&=& f(z). \\
\end{eqnarray*}

\bibliographystyle{alpha}
\bibliography{mybib}

\begin{thebibliography}{Rai06b}

\bibitem[Ber92]{Ber92}
B.\ Berndtsson.
\newblock Weighted estimates for $\bar\partial$ in domains in $\mathbb{C}$.
\newblock {\em Duke Math.\ J.}, 66(2):239--255, 1992.

\bibitem[BR09]{BoRa09}
A.\ Boggess and A.\ Raich.
\newblock A simplified calculation for the fundamental solution to the {H}eat
  {E}quation on the {H}eisenberg {G}roup.
\newblock {\em Proc.\ Amer.\ Math.\ Soc.}, 137(3):937--944, 2009.

\bibitem[Chr91]{Christ91}
M.\ Christ.
\newblock On the $\bar\partial$ equation in weighted ${L}^2$ norms in
  ${{\mathbb C}}^1$.
\newblock {\em J.\ Geom.\ Anal.}, 1(3):193--230, 1991.

\bibitem[FS91]{FoSi91}
J.E.\ Forn{\ae}ss and N.\ Sibony.
\newblock On ${L}^p$ estimates for $\overline\partial$.
\newblock In {\em Several {C}omplex {V}ariables and {C}omplex {G}eometry,
  {P}art 3 ({S}anta {C}ruz, {CA}, 1989)}, Proc.\ Sympos.\ Pure Math., 52, Part
  3, pages 129--163, Providence, R.I., 1991. American Mathematical Society.

\bibitem[Has94]{Has94}
F.~Haslinger.
\newblock Szeg{\"o} kernels for certain unbounded domains in ${{\mathbb
  C}}\sp2$. {T}ravaux de la {C}onf{\'e}rence {I}nternationale d{'}{A}nalyse
  {C}omplexe et du 7e {S}{\'e}minaire {R}oumano-{F}inlandais (1993).
\newblock {\em Rev.\ Roumaine Math.\ Pures Appl.}, 39:939--950, 1994.

\bibitem[Has95]{Has95}
F.~Haslinger.
\newblock Singularities of the {S}zeg{\"o} kernel for certain weakly
  pseudoconvex domains in ${C}\sp 2$.
\newblock {\em J.\ Funct.\ Anal.}, 129:406--427, 1995.

\bibitem[Has98]{Has98}
F.\ Haslinger.
\newblock {B}ergman and {H}ardy spaces on model domains.
\newblock {\em Illinois J.\ Math.}, 42:458--469, 1998.

\bibitem[Hul76]{Hul76}
A.\ Hulanicki.
\newblock The distribution of energy in the {B}rownian motion in the {G}aussian
  field and analytic hypoellipticity of certain subelliptic operators on the
  {H}eisenberg group.
\newblock {\em Studia Math.}, 56:165--173, 1976.

\bibitem[PR03]{PeRi03}
Marco~M.\ Peloso and Fulvio Ricci.
\newblock Analysis of the {K}ohn {L}aplacian on quadratic {CR} manifolds.
\newblock {\em J.\ Funct.\ Anal.}, 2003(2):321--355, 2003.

\bibitem[Rai]{Rai09h}
Andrew Raich.
\newblock Heat equations and the weighted $\dbar$-problem with decoupled
  weights.
\newblock {\em submitted}.
\newblock arXiv:0704.2768.

\bibitem[Rai06a]{Rai06h}
Andrew Raich.
\newblock Heat equations in ${{\mathbb R}}\times{{\mathbb C}}$.
\newblock {\em J.~Funct.\ Anal.}, 240(1):1--35, 2006.

\bibitem[Rai06b]{Rai06f}
Andrew Raich.
\newblock One-parameter families of perators in ${\mathbb{\C}}$.
\newblock {\em J.\ Geom.\ Anal.}, 16(2):353--374, 2006.

\bibitem[Rai07]{Rai07}
Andrew Raich.
\newblock Pointwise estimates of relative fundamental solutions for heat
  equations in ${{\mathbb R}}\times{{\mathbb C}}$.
\newblock {\em Math.\ Z.}, 256:193--220, 2007.

\bibitem[Ste93]{Ste93}
Elias~M.\ Stein.
\newblock {\em Harmonic Analysis: Real-Variable Methods, Orthogonality, and
  Oscillatory Integrals}.
\newblock Princeton Mathematical Series; 43. Princeton University Press,
  Princeton, New Jersey, 1993.

\bibitem[Tha93]{Tha93}
Sundaram Thangavelu.
\newblock {\em Lectures on Hermite and Laguerre Expansions}, volume~42 of {\em
  Mathematical Notes}.
\newblock Princeton University Press, Princeton, New Jersey, 1993.

\bibitem[YZ08]{YaZh08}
Qiaohua Yang and Fuliu Zhu.
\newblock The heat kernel on h-type groups.
\newblock {\em Proc.\ Amer.\ Math.\ Soc.}, 136(4):1457--1464, 2008.

\end{thebibliography}

\end{document}